\documentclass[11pt]{amsart}
\usepackage{mathrsfs}
\usepackage{amssymb,latexsym}
\usepackage{pstcol,pstricks,color}

 \numberwithin{equation}{section}

\newtheorem{theorem}{Theorem}[section]

\newtheorem{corollary}[theorem]{Corollary}
\newtheorem{definition}[theorem]{Definition}
\newtheorem{conjecture}[theorem]{Conjecture}
\newtheorem{remark}[theorem]{Remark}

\newtheorem{lemma}[theorem]{Lemma}

\newtheorem{example}[theorem]{Example}

\newtheorem{question}[theorem]{Question}

\def\qed{\hfill $\Box$}
\def\pf{\noindent {\it Proof.} }

\title{The Star of David Rule}

\begin{document}
\maketitle
\begin{center}
Yidong Sun

Department of Mathematics, Dalian Maritime University, 116026 Dalian, P.R. China\\[5pt]

{\it  sydmath@yahoo.com.cn}
\end{center}\vskip0.5cm

\subsection*{Abstract} In this note, a new concept
called {\em $SDR$-matrix} is proposed, which is an infinite lower
triangular matrix obeying the generalized rule of David star. Some
basic properties of $SDR$-matrices are discussed and two conjectures
on $SDR$-matrices are presented, one of which states that if a
matrix is a $SDR$-matrix, then so is its matrix inverse (if exists).

\medskip

{\bf Keywords}:  Narayana triangle, Pascal triangle, Lah triangle,
$SDR$-matrix.

\noindent {\sc 2000 Mathematics Subject Classification}: Primary
05A10; Secondary 15A09

\section{Introduction}

The \emph{Star of David rule} ~\cite{web}, originally stated by
Gould in 1972, is given by
\begin{equation*}
\binom{n}{k}\binom{n+1}{k-1}\binom{n+2}{k+1}=\binom{n}{k-1}\binom{n+1}{k+1}
\binom{n+2}{k},
\end{equation*}
for any $k$ and $n$, which implies that
\begin{equation*}
\binom{n}{k+1}\binom{n+1}{k}\binom{n+2}{k+2}=\binom{n}{k}\binom{n+1}{k+2}
\binom{n+2}{k+1}.
\end{equation*}
In 2003, the author observed in his Master dissertation \cite{Sun}
that if multiplying the above two identities and dividing by
$n(n+1)(n+2)$, one can arrive at
\begin{equation*}
N_{n,k+1}N_{n+1,k}N_{n+2,k+2}=N_{n,k}N_{n+1,k+2}N_{n+2,k+1},
\end{equation*}
where $N_{n,k}=\frac{1}{n}\binom{n}{k}\binom{n}{k-1}$ is the
Narayana number \cite[A001263]{sloane}.

In the summer of 2006, the author asked Mansour \cite{Mansour} for a
combinatorial proof of the above Narayana identity to be found.
Later, by Chen's bijective algorithm for trees \cite{Chen}, Li and
Mansour \cite{LiMan} provided a combinatorial proof of a general
identity
\begin{eqnarray*}
&&N_{n,k+m-1}N_{n+1,k+m-2}N_{n+2,k+m-3}\cdots N_{n+m-2,k+1}N_{n+m-1,k}N_{n+m,k+m}\qquad  \\
&&\qquad \qquad \qquad =N_{n,k}N_{n+1,k+m}N_{n+2,k+m-1}\cdots
N_{n+m-2,k+3}N_{n+m-1,k+2}N_{n+m,k+1}.
\end{eqnarray*}

This motivates the author to reconsider the Star of David rule and
to propose a new concept called {\em $SDR$-matrix} which obeys the
generalized rule of David star.

\begin{definition}\label{defi 1.1}
Let $\mathscr{A}=\Big(A_{n,k}\Big)_{n\geq k\geq 0}$ be an infinite
lower triangular matrix, for any given integer $m\geq 3$, if there
hold
\begin{eqnarray*}
\prod_{i=0}^rA_{n+i, k+r-i}\prod_{i=0}^{p-r-1}A_{n+p-i, k+r+i+1}=
\prod_{i=0}^{r}A_{n+p-i, k+p-r+i}\prod_{i=0}^{p-r-1}A_{n+i,
k+p-r-i-1},
\end{eqnarray*}
for all $2\leq p\leq m-1$ and $0\leq r\leq p-1$, then $\mathscr{A}$
is called an {\em $SDR$-matrix of order $m$}.
\end{definition}

\begin{figure}[h]
\setlength{\unitlength}{0.4mm}
\begin{center}
\begin{pspicture}(13,5.3)
\psset{xunit=25pt,yunit=25pt}\psgrid[subgriddiv=1,griddots=5,
gridlabels=4pt](0,0)(15,6)

\psline(1,4)(2,2)(3,3)(1,4)\psline(1,3)(3,2)(2,4)(1,3)

\pscircle*(1,4){0.06}\pscircle*(1,3){0.06}\pscircle*(1,2){0.06}
\pscircle*(2,4){0.06}\pscircle*(2,3){0.06}\pscircle*(2,2){0.06}
\pscircle*(3,4){0.06}\pscircle*(3,3){0.06}\pscircle*(3,2){0.06}

\psline[linewidth=1.2pt](5,4.5)(6,1.5)(8,3.5)(5,4.5)\psline[linewidth=1.2pt](5,2.5)(8,1.5)(7,4.5)(5,2.5)
\psline[linewidth=.5pt](5,3.5)(6,4.5)(7,1.5)(8,2.5)(5,3.5)

\pscircle*(5,1.5){0.06}\pscircle*(5,4.5){0.06}\pscircle*(5,3.5){0.06}\pscircle*(5,2.5){0.06}
\pscircle*(6,1.5){0.06}\pscircle*(6,4.5){0.06}\pscircle*(6,3.5){0.06}\pscircle*(6,2.5){0.06}
\pscircle*(7,1.5){0.06}\pscircle*(7,4.5){0.06}\pscircle*(7,3.5){0.06}\pscircle*(7,2.5){0.06}
\pscircle*(8,1.5){0.06}\pscircle*(8,4.5){0.06}\pscircle*(8,3.5){0.06}\pscircle*(8,2.5){0.06}

\psline[linewidth=1.2pt](10,5)(11,1)(14,4)(10,5)\psline[linewidth=1.2pt](14,1)(10,2)(13,5)(14,1)
\psline[linewidth=.5pt](10,4)(11,5)(14,3)(12,1)(10,4)\psline[linewidth=.5pt](10,3)(12,5)(14,2)(13,1)(10,3)

\pscircle*(10,5){0.06}\pscircle*(10,5){0.06}\pscircle*(10,4){0.06}\pscircle*(10,3){0.06}\pscircle*(10,2){0.06}\pscircle*(10,1){0.06}
\pscircle*(11,5){0.06}\pscircle*(11,5){0.06}\pscircle*(11,4){0.06}\pscircle*(11,3){0.06}\pscircle*(11,2){0.06}\pscircle*(11,1){0.06}
\pscircle*(12,5){0.06}\pscircle*(12,5){0.06}\pscircle*(12,4){0.06}\pscircle*(12,3){0.06}\pscircle*(12,2){0.06}\pscircle*(12,1){0.06}
\pscircle*(13,5){0.06}\pscircle*(13,5){0.06}\pscircle*(13,4){0.06}\pscircle*(13,3){0.06}\pscircle*(13,2){0.06}\pscircle*(13,1){0.06}
\pscircle*(14,5){0.06}\pscircle*(14,5){0.06}\pscircle*(14,4){0.06}\pscircle*(14,3){0.06}\pscircle*(14,2){0.06}\pscircle*(14,1){0.06}

\put(1.3,.3){$p=2$}\put(5.25,.3){$p=3$}\put(10.1,.3){$p=4$}

\end{pspicture}
\caption{The case $m=5$. }\label{fDD1}
\end{center}
\end{figure}
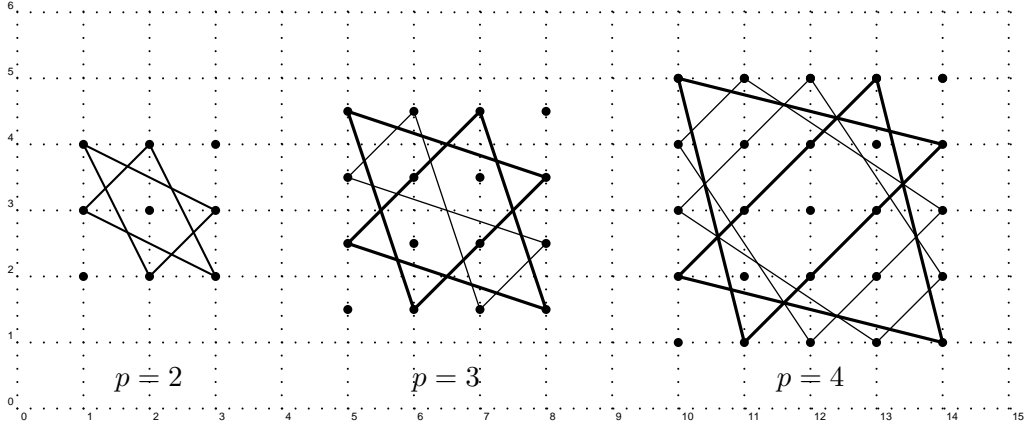

In order to give a more intuitive view on the definition, we present
a pictorial description of the generalized rule for the case $m=5$.
See Figure \ref{fDD1}.

Let $SDR_{m}$ denote the set of $SDR$-matrices of order $m$ and
$SDR_{\infty}$ be the set of $SDR$-matrices $\mathscr{A}$ of order
$\infty$, that is $\mathscr{A}\in SDR_{m}$ for any $m\geq 3$. By our
notation, it is obvious that the Pascal triangle
$\mathscr{P}=\Big(\binom{n}{k}\Big)_{n\geq k\geq 0}$ and the
Narayana triangle $\mathscr{N}=\Big(N_{n+1,k+1}\Big)_{n\geq k\geq
0}$ are $SDR$-matrices of order $3$. In fact, both of them will be
proved to be $SDR$-matrices of order $\infty$.
\begin{eqnarray*}
\begin{array}{cc}
\mathscr{P}=\left(
\begin{array}{rrrrrr}
 1&    &    &    &     &   \\
 1&   1&    &    &     &   \\
 1&   2&   1&    &     &   \\
 1&   3&   3&   1&     &   \\
 1&   4&   6&   4&    1&   \\
 1&   5&  10&  10&    5& 1 \\
  &    &\cdots &    &     &   \\
\end{array}\right), & \hskip0.2cm
\mathscr{N}=\left(
\begin{array}{rrrrrr}
 1&    &    &    &     &   \\
 1&   1&    &    &     &   \\
 1&   3&   1&    &     &   \\
 1&   6&   6&   1&     &   \\
 1&  10&  20&  10&    1&   \\
 1&  15&  50&  50&   15& 1 \\
  &    &\cdots &    &     &   \\
\end{array}\right).
\end{array}
\end{eqnarray*}

In this paper, we will discuss some basic properties of the sets
$SDR_m$ and propose two conjectures on $SDR_m$ for $3\leq m\leq
\infty$ in the next section. We also give some comments on relations
between $SDR$-matrices and Riordan arrays in Section $3$.

\section{The basic properties of $SDR$-matrices}
For any infinite lower triangular matrices
$\mathscr{A}=\Big(A_{n,k}\Big)_{n\geq k\geq 0}$ and
$\mathscr{B}=\Big(B_{n,k}\Big)_{n\geq k\geq 0}$, define
$\mathscr{A}\circ\mathscr{B}=\Big(A_{n,k}B_{n,k}\Big)_{n\geq k\geq
0}$ to be the Hadamard product of $\mathscr{A}$ and $\mathscr{B}$,
denote by $\mathscr{A}^{\circ j}$ the $j$-th Hadamard power of
$\mathscr{A}$; If $A_{n,k}\neq 0$ for $n\geq k\geq 0$, then define
$\mathscr{A}^{\circ(-1)}=\Big(A_{n,k}^{-1}\Big)_{n\geq k\geq 0}$ to
be the Hadamard inverse of $\mathscr{A}$.



From Definition \ref{defi 1.1}, one can easily derive the following
three lemmas.

\begin{lemma}\label{lemma 2.1}
For any $\mathscr{A}\in SDR_{m}$, $\mathscr{B}\in SDR_{m+i}$ with
$i\geq 0$, there hold $\mathscr{A}\circ\mathscr{B}\in SDR_{m}$, and
$\mathscr{A}^{\circ(-1)}\in SDR_{m}$ if it exists.
\end{lemma}



\begin{lemma}\label{lemma 2.2}
For any $\mathscr{A}=\Big(A_{n,k}\Big)_{n\geq k\geq 0}\in SDR_{m}$,
then $\Big(A_{n+i,k+j}\Big)_{n\geq k\geq 0}\in SDR_{m}$ for fixed
$i, j\geq 0$.
\end{lemma}

\begin{lemma}\label{lemma 2.3}
Given any sequence $(a_n)_{n\geq 0}$, let $A_{n,k}=a_n$,
$B_{n,k}=a_k$ and $C_{n,k}=a_{n-k}$ for $n\geq k\geq 0$, then
$\Big(A_{n,k}\Big)_{n\geq k\geq 0}, \Big(B_{n,k}\Big)_{n\geq k\geq
0}, \Big(C_{n,k}\Big)_{n\geq k\geq 0}\in SDR_{\infty}$.
\end{lemma}

\begin{example}{\rm Let $a_n=n!$ for $n\geq 0$, then we have
\begin{eqnarray*}
\mathscr{P} &=& \Big(n!\Big)_{n\geq k\geq 0}\circ
\Big(k!\Big)_{n\geq k\geq 0}^{\circ(-1)}\circ
\Big((n-k)!\Big)_{n\geq k\geq 0}^{\circ(-1)}, \\
\mathscr{N} &=& \Big(\frac{1}{k+1}\Big)_{n\geq k\geq 0}\circ
\mathscr{P}\circ \Big(\binom{n+1}{k}\Big)_{n\geq k\geq 0}, \\
\mathscr{L} &=& \Big((n+1)!\Big)_{n\geq k\geq
0}\circ\mathscr{P}\circ \Big((k+1)!\Big)_{n\geq k\geq
0}^{\circ(-1)},
\end{eqnarray*}
which, by Lemmas \ref{lemma 2.1}-\ref{lemma 2.3}, produce that the
Pascal triangle $\mathscr{P}$, the Narayana triangle $\mathscr{N}$
and the Lah triangle $\mathscr{L}$ belong to $SDR_{\infty}$, where
$(\mathscr{L})_{n,k}=\binom{n}{k}\frac{(n+1)!}{(k+1)!}$ is the Lah
number \cite{comtet}. }
\end{example}

\begin{theorem}\label{theo 0}
For any sequences $(a_n)_{n\geq 0}$, $(b_n)_{n\geq 0}$ and
$(c_n)_{n\geq 0}$ such that $b_0=1$, $a_n\neq 0$ and $c_n\neq 0$ for
$n\geq 0$, let $\mathscr{A}=\Big(a_kb_{n-k}c_n\Big)_{n\geq k\geq
0}$, then $\mathscr{A}^{-1}\in SDR_{\infty}$.
\end{theorem}
\pf By Lemmas \ref{lemma 2.1} and \ref{lemma 2.3}, we have
$\mathscr{A}\in SDR_{\infty}$. It is not difficult to derive the
matrix inverse $\mathscr{A}^{-1}$ of $\mathscr{A}$ with the generic
entries
\begin{eqnarray*}
\Big(\mathscr{A}^{-1}\Big)_{n,k}&=&a_n^{-1}B_{n-k}c_{k}^{-1},
\end{eqnarray*}
where $B_n$ with $B_0=1$ are given by
\begin{eqnarray}\label{eqn 2.0}
B_n&=&\displaystyle\sum_{j=1}^n(-1)^j\sum_{i_1+i_2+\cdots+i_j=n,
i_1, \dots, i_j\geq 1}b_{i_1}b_{i_2}\cdots b_{i_j},  \ (n\geq 1).
\end{eqnarray}
Hence, by Lemmas \ref{lemma 2.1} and \ref{lemma 2.3}, one can deduce
that
\begin{eqnarray*}
\mathscr{A}^{-1} &=& \Big(a_n^{-1}\Big)_{n\geq k\geq
0}\circ\Big(B_{n-k}\Big)_{n\geq k\geq 0}\circ
\Big(c_k^{-1}\Big)_{n\geq k\geq 0}\in SDR_{\infty},
\end{eqnarray*}
as desired. \qed\vskip0.2cm

Specially, when $c_n:=1$ or $a_n:=\frac{a_n}{n!}$,
$b_n:=\frac{b_n}{n!}$, $c_n:=n!$, both
$\mathscr{B}=\Big(a_kb_{n-k}\Big)_{n\geq k\geq 0}$ and
$\mathscr{C}=\Big(\binom{n}{k}a_kb_{n-k}\Big)_{n\geq k\geq 0}$ are
in $SDR_{\infty}$, then so $\mathscr{B}^{-1}$ and
$\mathscr{C}^{-1}$. More precisely, let
$a_n^{-1}=b_n^{-1}=c_n=n!(n+1)!$ for $n\geq 0$, note that the
Narayana triangle $\mathscr{N}\in SDR_{\infty}$ and
\begin{eqnarray*}
N_{n+1,k+1}=\frac{1}{n+1}\binom{n+1}{k+1}\binom{n+1}{k}=\frac{n!(n+1)!}{k!(k+1)!(n-k)!(n-k+1)!}.
\end{eqnarray*}
Then one has $\mathscr{N}^{-1}\in SDR_{\infty}$ by Theorem \ref{theo
0}.

Theorem \ref{theo 0} suggests the following conjecture.
\begin{conjecture}
For any $\mathscr{A}\in SDR_{m}$, if the inverse $\mathscr{A}^{-1}$
of $\mathscr{A}$ exists, then $\mathscr{A}^{-1}\in SDR_{m}$.
\end{conjecture}

\begin{theorem}\label{theo 1}
For any sequences $(a_n)_{n\geq 0}$, $(b_n)_{n\geq 0}$ with $b_0=1$
and $a_n\neq 0$ for $n\geq 0$, let
$\mathscr{A}=\Big(a_nb_{n-k}a_k^{-1}\Big)_{n\geq k\geq 0}$, then the
matrix power $\mathscr{A}^{j}\in SDR_{\infty}$ for any integer $j$.
\end{theorem}
\pf By Lemmas \ref{lemma 2.1} and \ref{lemma 2.3}, we have
$\mathscr{A}\in SDR_{\infty}$. Note that it is trivially true for
$j=1$ and $j=0$ (where $\mathscr{A}^0$ is the identity matrix by
convention). It is easy to obtain the $(n,k)$-entries of
$\mathscr{A}^{j}$ for $j\geq 2$,
\begin{eqnarray*}
\Big(\mathscr{A}^{j}\Big)_{n,k}&=&\sum_{k\leq k_{j-1}\leq \cdots\leq
k_1\leq
n}\mathscr{A}_{n,k_1}\mathscr{A}_{k_1,k_2}\cdots\mathscr{A}_{k_{j-2},k_{j-1}}\mathscr{A}_{k_{j-1},k}
\\
&=&a_{n}C_{n-k}a_{k}^{-1},
\end{eqnarray*}
where $C_n$ with $C_0=1$ is given by
$C_n=\sum_{i_1+i_2+\cdots+i_j=n, i_1,\dots,i_j\geq 0}b_{i_1}
b_{i_2}\cdots b_{i_{j}}$ for $n\geq 1$.

By Lemmas \ref{lemma 2.1} and \ref{lemma 2.3}, one can deduce that
\begin{eqnarray*}
\mathscr{A}^{j} &=& \Big(a_n\Big)_{n\geq k\geq
0}\circ\Big(C_{n-k}\Big)_{n\geq k\geq 0}\circ
\Big(a_k^{-1}\Big)_{n\geq k\geq 0}\in SDR_{\infty}.
\end{eqnarray*}
By Theorem \ref{theo 0} and its proof, we have $\mathscr{A}^{-1}\in
SDR_{\infty}$ and
$\big(\mathscr{A}^{-1}\big)_{n,k}=a_nB_{n-k}a_{k}^{-1}$, where $B_n$
is given by (\ref{eqn 2.0}). Note that $\mathscr{A}^{-1}$ has the
form as required in Theorem \ref{theo 1}, so by the former part of
this proof, we have $\mathscr{A}^{-j}\in SDR_{\infty}$ for $j\geq
1$. Hence we are done. \qed \vskip0.2cm

Let $a_n=b_n=n!$, $a_n=b_n=n!(n+1)!$ or $a_n=n!(n+1)!$ and
$b_n^{-1}=n!$ for $n\geq 0$ in Theorem \ref{theo 1}, one has
\begin{corollary}
For $\mathscr{P}$, $\mathscr{N}$ and $\mathscr{L}$, then
$\mathscr{P}^j, \mathscr{N}^j, \mathscr{L}^j\in SDR_{\infty}$ for
any integer $j$.
\end{corollary}

\begin{remark}
{\rm In general, for $\mathscr{A}, \mathscr{B}\in SDR_{m}$, their
matrix product $\mathscr{A}\mathscr{B}$ is possibly not in
$SDR_{m}$. For example, $\mathscr{P},\mathscr{N}\in SDR_{3}$, but
\begin{eqnarray*}
\mathscr{P}\mathscr{N}&=&\left(
\begin{array}{rrrrrr}
 1&    &    &    &     &   \\
 2&   1&    &    &     &   \\
 4&   5&   1&    &     &   \\
 8&  18&   9&   1&     &   \\
16&  56&  50&  14&    1&   \\
32& 160& 220& 110&   20& 1 \\
  &    &\cdots &    &     &   \\
\end{array}\right)
\notin SDR_{3}.
\end{eqnarray*} }
\end{remark}

\begin{theorem}\label{theo 2}
For any $\mathscr{A}=\Big(A_{n,k}\Big)_{n\geq k\geq 0}$ with
$A_{n,k}\neq 0$ for $n\geq k\geq 0$, then $\mathscr{A}\in SDR_{m+1}$
if and only if  $\mathscr{A}\in SDR_{m}$.
\end{theorem}
\pf Note that $SDR_{m+1}\subset SDR_{m}$, so the necessity is clear.
It only needs to prove the sufficient condition. For the symmetry,
it suffices to verify
\begin{eqnarray*}
\prod_{i=0}^rA_{n+i, k+r-i}\prod_{i=0}^{m-r}A_{n+m-i+1, k+r+i+1}=
\prod_{i=0}^{r}A_{n+m-i+1, k+m-r+i+1}\prod_{i=0}^{m-r}A_{n+i,
k+m-r-i},
\end{eqnarray*}
for $0\leq r\leq [m/2]-1$. We just take the case $r=0$ for example,
others can be done similarly. It is trivial when $A_{n, k+m}=A_{n+1,
k+m+1}=0$. So we assume that $A_{n, k+m}\neq 0, A_{n+1, k+m+1}\neq
0$, then all $A_{n+i, k+j}$ to be considered, except for $A_{n,
k+m+1}$, must not be zero. By Definition \ref{defi 1.1}, we have
\begin{eqnarray}\label{eqn 2.1}
\lefteqn{A_{n+m-i, k+i}A_{n+m-i-1, k+i+1}A_{n+m-i+1, k+i+2} } \nonumber\\
&=&A_{n+m-i+1, k+i+1}A_{n+m-i, k+i+2}A_{n+m-i-1, k+i}, \hskip0.2cm
(0\leq i\leq m-1).
\end{eqnarray}
\begin{eqnarray}\label{eqn 2.2}
A_{n+m+1, k+m+1}\prod_{i=0}^{m-1}A_{n+i, k+m-i}=A_{n+1,
k+1}\prod_{i=0}^{m-1}A_{n+m-i+1, k+i+2}.
\end{eqnarray}
\begin{eqnarray}\label{eqn 2.3}
A_{n+1, k+1}\prod_{i=0}^{m-1}A_{n+m-i, k+i+1}=A_{n+m,
k+m}\prod_{i=0}^{m-1}A_{n+i+1, k+m-i-1}.
\end{eqnarray}
\begin{eqnarray}\label{eqn 2.4}
A_{n+m, k+m}\prod_{i=0}^{m-1}A_{n+i, k+m-i-1}=A_{n,
k}\prod_{i=0}^{m-1}A_{n+m-i, k+i+1}.
\end{eqnarray}
Multiplying (\ref{eqn 2.1})$-$(\ref{eqn 2.4}) together, after
cancellation, one can get
\begin{eqnarray*}
A_{n, k}\prod_{i=0}^{m}A_{n+m-i+1, k+i+1}=A_{n+m+1,
k+m+1}\prod_{i=0}^{m}A_{n+i, k+m-i},
\end{eqnarray*}
which confirms the case $r=0$. \qed \vskip0.2cm

\begin{remark}
{\rm The condition $A_{n,k}\neq 0$ for $n\geq k\geq 0$ in Theorem
\ref{theo 2} is necessary. The following example verifies this
claim.
\begin{eqnarray*}
\Big(\binom{\frac{n+k}{2}}{\frac{n-k}{2}}\Big)_{n\geq k\geq
0}=\left(
\begin{array}{cccccc}
 1&    &    &    &     &   \\
 0&   1&    &    &     &   \\
 1&   0&   1&    &     &   \\
 0&   2&   0&   1&     &   \\
 1&   0&   3&   0&    1&   \\
 0&   3&   0&   4&    0& 1 \\
  &    &\cdots &    &     &   \\
\end{array}\right)
\in SDR_{3}, \mbox{but\ not\ in}\ SDR_{4}.
\end{eqnarray*} }
\end{remark}

Recall that the Narayana number $\mathscr{N}_{n+1,k+1}$ can be
represented as
\begin{eqnarray*}
\mathscr{N}_{n+1,k+1}=\frac{1}{n+1}\binom{n+1}{k+1}\binom{n+1}{k}=
\det\left(
\begin{array}{cc}
\binom{n}{k} & \binom{n}{k+1} \\[5pt]
\binom{n+1}{k} & \binom{n+1}{k+1}
\end{array} \right),
\end{eqnarray*}
so we can come up with the following definition.
\begin{definition}\label{defi 2.1}
Let $\mathscr{A}=\Big(A_{n,k}\Big)_{n\geq k\geq 0}$ be an infinite
lower triangular matrix, for any integer $j\geq 1$, define
$\mathscr{A}_{[j]}=\Big(A_{n,k}^{[j]}\Big)_{n\geq k\geq 0}$, where
\begin{eqnarray*}
A_{n,k}^{[j]}=\det\left(
\begin{array}{ccc}
A_{n,k}     &  \cdots  & A_{n,k+j-1} \\[5pt]
 \vdots     &  \cdots  & \vdots       \\[5pt]
A_{n+j-1,k} &  \cdots  & A_{n+j-1,k+j-1}
\end{array} \right).
\end{eqnarray*}
\end{definition}

\begin{theorem}\label{theo 3}
For any sequences $(a_n)_{n\geq 0}$, $(b_n)_{n\geq 0}$ and
$(c_n)_{n\geq 0}$ such that $b_0=1$, $a_n\neq 0$ and $c_n\neq 0$ for
$n\geq 0$, let $\mathscr{A}=\Big(a_kb_{n-k}c_n\Big)_{n\geq k\geq
0}$, then $\mathscr{A}_{[j]}\in SDR_{\infty}$ for any integer $j\geq
1$.
\end{theorem}
\pf By Lemmas \ref{lemma 2.1} and \ref{lemma 2.3}, we have
$\mathscr{A}\in SDR_{\infty}$. It is easy to derive the determinant
\begin{eqnarray*}
\det\left(
\begin{array}{ccc}
a_kb_{n-k}c_n            &  \cdots  & a_{k+j-1}b_{n-k-j+1}c_n \\[5pt]
 \vdots                  &  \cdots  & \vdots          \\[5pt]
a_kb_{n-k+j-1}c_{n+j-1}  &  \cdots  & a_{k+j-1}b_{n-k}c_{n+j-1}
\end{array}\right)
=B_{n-k}\prod_{i=0}^{j-1}a_{k+i}c_{n+i},
\end{eqnarray*}
where $B_n$ with $B_0=1$ are given by
\begin{eqnarray*}
B_n=\det\left(
\begin{array}{ccc}
b_{n}      &  \cdots  &  b_{n-j+1} \\[5pt]
 \vdots    &  \cdots  &  \vdots    \\[5pt]
b_{n+j-1}  &  \cdots  &  b_{n}
\end{array}\right).
\end{eqnarray*}
Hence, by Lemmas \ref{lemma 2.1} and \ref{lemma 2.3}, one can deduce
that
\begin{eqnarray*}
\mathscr{A}_{[j]} &=& \Big(\prod_{i=0}^{j-1}a_{k+i}\Big)_{n\geq
k\geq 0}\circ\Big(B_{n-k}\Big)_{n\geq k\geq 0}\circ
\Big(\prod_{i=0}^{j-1}c_{n+i}\Big)_{n\geq k\geq 0}\in SDR_{\infty},
\end{eqnarray*}
as desired. \qed\vskip0.2cm

Let $a_n^{-1}=b_n^{-1}=c_n=n!$, $a_n^{-1}=b_n^{-1}=c_n=n!(n+1)!$ or
$a_n^{-1}=c_n=n!(n+1)!$ and $b_n^{-1}=n!$ for $n\geq 0$ in Theorem
\ref{theo 3}, one has
\begin{corollary}
For $\mathscr{P}$, $\mathscr{N}$ and $\mathscr{L}$, then
$\mathscr{P}_{[j]}, \mathscr{N}_{[j]}, \mathscr{L}_{[j]}\in
SDR_{\infty}$ for any integer $j\geq 1$.
\end{corollary}

Theorem \ref{theo 3} suggests the following conjecture.
\begin{conjecture}
If $\mathscr{A}\in SDR_{\infty}$, then $\mathscr{A}_{[j]}\in
SDR_{\infty}$ for any integer $j\geq 1$.
\end{conjecture}

\begin{remark}
{\rm The conjecture on $SDR_{m}$ is generally not true for $3\leq m<
\infty$. For example, let $\mathscr{A}=\Big(A_{n,k}\Big)_{n\geq
k\geq 0}$ with $A_{n,k}=\binom{\frac{n+k}{2}}{\frac{n-k}{2}}$, then
we have $\mathscr{A}\in SDR_{3}$, but
\begin{eqnarray*}
\begin{array}{rr}
\mathscr{A}_{[2]}=\left(
\begin{array}{rrrrr}
             1&                &                &    &          \\
            -1&               1&                &    &          \\
         \bf 2&     -\textrm{2}&               1&    &          \\
   -\textrm{2}&               6&           \bf-3&   1&          \\
             3&           \bf-9&     \textrm{12}&  -4&    1     \\
              &                &         \cdots &    &
\end{array}\right)
\notin SDR_{3}, & \mathscr{A}_{[3]}=\left(
\begin{array}{rrrrr}
             1&                &                &    &          \\
             0&               1&                &    &          \\
             2&               0&               1&    &          \\
             0&              15&               0&   1&          \\
             9&               0&              36&   0&    1    \\
              &                &         \cdots &    &
\end{array}\right)
\in SDR_{3}.
\end{array}
\end{eqnarray*} }
\end{remark}

\vskip0.5cm

\section{Further Comments}

We will present some further comments on the connections between
$SDR$-matrices and Riordan arrays. The concept of Riordan array
introduced by Shapiro et al \cite{SGWW}, plays a particularly
important role in studying combinatorial identities or sums and also
is a powerful tool in study of many counting problems
\cite{MRSV,MSV,MV}. For examples, Sprugnoli \cite{MSV,Sp1,Sp2}
investigated Riordan arrays related to binomial coefficients,
colored walks, Stirling numbers and Abel-Gould identities.

To define a Riordan array we need two analytic functions,
$d(t)=d_0+d_1t+d_2t^2+\cdots$ and $h(t)=h_1t+h_2t^2+\cdots$. A {\em
Riordan array} is an infinite lower triangular array
$\{d_{n,k}\}_{n,k\in \mathbb{N}}$, defined by a pair of formal power
series $(d(t),h(t))$, with the generic element $d_{n,k}$ satisfying
\begin{eqnarray*}
d_{n,k}&=&[t^n]d(t)(h(t))^k, \ \ \ (n,\ k\geq 0).
\end{eqnarray*}
Assume that $d_0\neq 0\neq h_1$, then $(d(t), h(t))$ is an element
of the {\em Riordan group} \cite{SGWW}, under the group
multiplication rule:
\begin{eqnarray*}
(d(t),h(t))(g(t),f(t))=(d(t)g(h(t)),f(h(t))).
\end{eqnarray*}
This indicates that the identity is $I=(1,t)$, the usual matrix
identity, and that
\begin{eqnarray*}
(d(t),h(t))^{-1}=(\frac{1}{d(\overline{h}(t))},\overline{h}(t)),
\end{eqnarray*}
where $\overline{h}(t)$ is the compositional inverse of $h(t)$,
i.e., $\overline{h}(h(t))=h(\overline{h}(t))=t$.

By our notation, we have
\begin{eqnarray*}
\mathscr{P}&=&(\frac{1}{1-t}, \frac{t}{1-t})\in SDR_{\infty},\\
\mathscr{P}^j&=&(\frac{1}{1-jt}, \frac{t}{1-jt})\in SDR_{\infty},\\
\Big(\binom{\frac{n+k}{2}}{\frac{n-k}{2}}\Big)_{n\geq k\geq 0}
  &=& (\frac{1}{1-t^2}, \frac{t}{1-t^2})\in SDR_{3},\\
(\frac{1}{1-t^2}, \frac{t}{1-t^2})^{-1}&=& (\frac{1-\sqrt{1-4t^2}}{2t^2}, \frac{1-\sqrt{1-4t^2}}{2t})\in SDR_{3},\\
\Big(d_{n-k}\Big)_{n\geq k\geq 0} &=& (d(t), t)\in SDR_{\infty}.
\end{eqnarray*}
Hence, it is natural to ask the following question.
\begin{question}
Given a formal power series $d(t)$, what conditions $h(t)$ should
satisfy, such that $(d(t), h(t))$ forms an $SDR$-matrix.
\end{question}

\vskip0.5cm

\section*{Acknowledgements} The author is grateful to the anonymous referees
for the helpful suggestions and comments. This work was supported by
The National Science Foundation of China.


\end{document}